\documentclass{amsart}
\usepackage{amssymb,latexsym,amsmath}
\usepackage[mathscr]{eucal}
\usepackage{amsthm}
\usepackage{verbatim}
 \newtheorem*{vartheorem*}{ }

\def\sideremark#1{\ifvmode\leavevmode\fi\vadjust{\vbox to0pt{\vss
\hbox to 0pt{\hskip\hsize\hskip1em
\vbox{\hsize2cm\tiny\raggedright\pretolerance10000
\noindent#1\hfill}\hss}\vbox to8pt{\vfil}\vss}}}

\newtheorem*{structure}{Structure Theorem}
\newtheorem*{theorem*}{Theorem}
\newtheorem*{corollary*}{Corollary}
\newtheorem{theorem}{Theorem}[section]
\newtheorem{corollary}[theorem]{Corollary}
\newtheorem{lemma}[theorem]{Lemma}
\newtheorem{proposition}[theorem]{Proposition}
\newtheorem{remark}[theorem]{Remark}
\newtheorem*{remark*}{Remark}
\newtheorem{definition}[theorem]{Definition}
\newtheorem{example}[theorem]{Example}
\newtheorem{question}[theorem]{Question}
\newcommand{\be}{\begin{equation}\label}
\newcommand{\ee}{\end{equation}}
\newcommand{\bq}{\begin{equation*}}
\newcommand{\eq}{\end{equation*}}
\newcommand{\ba}{\begin{align*}}
\newcommand{\ea}{\end{align*}}
\newcommand{\bp}{\begin{proof}}
\newcommand{\ep}{\end{proof}}
\newcommand{\bL}{\begin{lemma}\label}
\newcommand{\eL}{\end{lemma}}
\newcommand{\bP}{\begin{proposition}\label}
\newcommand{\eP}{\end{proposition}}
\newcommand{\bC}{\begin{corollary}\label}
\newcommand{\eC}{\end{corollary}}
\newcommand{\bT}{\begin{theorem}\label}
\newcommand{\eT}{\end{theorem}}
\newcommand{\bTT}{\begin{theorem*}\label}
\newcommand{\eTT}{\end{theorem*}}
\newcommand{\bR}{\begin{remark}\label}
\newcommand{\eR}{\end{remark}}
\newcommand{\bD}{\begin{definition}\label}
\newcommand{\eD}{\end{definition}}

\newcommand{\bE}{\begin{example}\label}
\newcommand{\eE}{\end{example}}
\newcommand{\bQ}{\begin{question}\label}
\newcommand{\eQ}{\end{question}}

\DeclareMathOperator{\diag}{diag}

\begin{document}
\title{SUBIDEALS OF OPERATORS II}
\author{Sasmita Patnaik}
\address{University of Cincinnati\\
          Department of Mathematics\\
          Cincinnati, OH, 45221-0025\\
          USA}
\email{sasmita\_19@yahoo.co.in}
\author{Gary Weiss}
\address{University of Cincinnati\\
          Department of Mathematics\\
          Cincinnati, OH, 45221-0025\\
          USA}
\email{gary.weiss@math.uc.edu}
\date{\today}

\begin{abstract}
A subideal (also called a $J$-ideal) is an ideal of a $B(H)$-ideal $J$.
This paper is the sequel to \textit{Subideals of Operators} where a complete characterization of principal and then finitely generated $J$-ideals were obtained by first generalizing the 1983 work of Fong and Radjavi who determined which principal $K(H)$-ideals are also $B(H)$-ideals. 
Here we determine which countably generated $J$-ideals are $B(H)$-ideals, and in the absence of the continuum hypothesis which $J$-ideals with generating sets of cardinality less than the continuum are $B(H)$-ideals. 
These and some other results herein are based on the dimension of a related quotient space. 
We use this to characterize these $J$-ideals and settle additional questions about subideals. 
A key property in our investigation turned out to be $J$-softness of a $B(H)$-ideal $I$ inside $J$, that is, $IJ = I$, a generalization of a recent notion of softness of $B(H)$-ideals introduced by  Kaftal-Weiss  and earlier exploited for Banach spaces by Mityagin and Pietsch.  
 
\end{abstract}

\maketitle

\section{Introduction}

In  \textit{Subideals of Operators} \cite{PW11} we found three types of principal and finitely generated subideals (i.e., $J$-ideals): linear, real-linear and nonlinear subideals. Such types also carry over to general $J$-ideals. 
The linear $K(H)$-ideals, being the traditional ones, were studied in 1983 by Fong and Radjavi \cite{FR83}.  
They found principal linear $K(H)$-ideals that are not $B(H)$-ideals. 
Herein we take all $J$-ideals $\mathcal I$ to be linear, but as proved in \cite{PW11}, we expect  here also that most of the results and methods apply to the two other types of subideals (real-linear and nonlinear). 
Also $H$, as in \cite{PW11}, denotes a separable infinite-dimensional complex Hilbert space. 
One of our main contributions in \cite{PW11} was to use a modern framework for $B(H)$-ideals to generalize \cite[Theorem 2]{FR83}. 
We generalized their result on principal $K(H)$-ideals to all principal $J$-ideals  by proving that a principal and then a finitely generated $J$-ideal $(\mathscr S)_J$ generated by the finite set $\mathscr S \subset J$ is also a $B(H)$-ideal if and only if $(\mathscr S)$ is $J$-soft, that is, $(\mathscr S)= (\mathscr S)J$ where $(\mathscr S)$ is the $B(H)$-ideal generated by $\mathscr S$. 
Then we used this to characterize the structure of $(\mathscr S)_J$. $J$-softness is a generalization of a recent notion of $K(H)$-softness of $B(H)$-ideals introduced by Kaftal-Weiss and earlier exploited for Banach spaces by Mityagin and Pietsch (see \cite[Remark 2.6]{PW11}, \cite{bM67}, \cite{aP80}).

Here we further develop the subject by investigating $J$-ideals $\mathcal I = (\mathscr S)_J$ generated by arbitrary sets $\mathscr S$ of varying cardinality, their algebraic structure and when they are $B(H)$-ideals. 
To add perspective, the reader should keep in mind that all nonzero $J$-ideals have cardinality and Hamel dimension precisely  cardinality $\mathfrak c$ of the continuum (Remark \ref{R:2}), but questions on the cardinalities of their possible generating sets is another matter (Section 6, Questions 1-2), and this we shall see impacts questions on structure. 
After investigating the cases when $\mathscr S$ is countable or of cardinality less than $\mathfrak c$ (absent the continuum hypothesis CH), we then consider general $J$-ideals $\mathcal I$ and questions on the possible cardinalities of their generating
sets, observing that $\mathcal I$ is always a generating set for itself but may have generating sets of cardinality less than its
cardinality $\mathfrak c$. 
When they do has special implications.  

We show $(\mathscr S)_J$ is a $B(H)$-ideal if and only if $(\mathscr S)$ is $J$-soft for those $(\mathscr S)_J$ generated as a $J$-ideal by countable sets and then when  generated by sets of cardinality strictly less than cardinality $\mathfrak c$  (Theorem \ref{T:1}).
This will follow from sufficiency of the codimension condition on $ (\mathscr S)_J^0$ in $ (\mathscr S)_J$: $ (\mathscr S)_J / (\mathscr S)_J^0$ has Hamel dimension strictly less than $\mathfrak c$. (For $ (\mathscr S)_J^0$ see Definition \ref{D:d2}.)
We then investigate general $J$-ideals  to provide an example where softness fails   for a $J$-ideal $\mathcal I$ with codimension of $\mathcal I^0$ in $\mathcal I$ equal to $\mathfrak c$ (Example \ref{E:11}), thereby showing that $J$-ideals that are also $B(H)$-ideals need not be $J$-soft and as a consequence cannot be generated in $J$ by sets of cardinality less than $\mathfrak c$. 
We also answer several questions on $J$-ideals posed in \cite[Sections 6-7]{PW11}, provide some additional results and pose new questions.

\vspace{.4cm}
\newpage

In summary the main theorems here are:\\

\noindent For $\mathcal I^0 := \text{span}\{\mathcal IJ + J\mathcal I\} + J(\mathcal I)J$ (see Definition \ref{D:d2}) where $(\mathcal I)$ is the $B(H)$-ideal generated by $J$-ideal $\mathcal I$,
\begin{theorem*}(Theorem \ref{T:1})
The $J$-ideal $(\mathscr S)_J$ generated by a set $\mathscr S$ of cardinality strictly less than $\mathfrak c$ is a $B(H)$-ideal if and only if  the $B(H)$-ideal $(\mathscr S)$ is $J$-soft (i.e., $(\mathscr S) = (\mathscr S)J$). Moreover, for a $J$-ideal $\mathcal I$ with Hamel dimension of $\mathcal I/\mathcal I^0$ strictly less than $\mathfrak c$, $\mathcal I$ is a $B(H)$-ideal if and only if $(\mathcal I)$ is $J$-soft.\\
\end{theorem*}

\begin{structure}(\text{Theorem \ref{T:3}})
For $(\mathscr S)_J$ when $|\mathscr S| < \mathfrak c$, \\
The algebraic structure of the $J$-ideal $(\mathscr S)_J$ generated by a set $\mathscr S$ is given by
\begin{center}
$(\mathscr S)_J = \text{span}\{\mathscr S + J\mathscr S + \mathscr SJ\} + J(\mathscr S)J$
\end{center}
$J(\mathscr S)J $ is a $B(H)$-ideal, $\text{span}\{J\mathscr S + \mathscr SJ\} + J(\mathscr S)J$ is a $J$-ideal, and $J(\mathscr S)J \subset \text{span}\{J\mathscr S + \mathscr SJ\} + J(\mathscr S)J \subset (\mathscr S)_J$.
This inclusion collapses to $J(\mathscr S)J = (\mathscr S)_J$ if and only if  $(\mathscr S)$ is $J$-soft.
\end{structure}

\begin{remark*}
Although the methods in \cite{FR83} are quite a bit more analytic, we found here and in \cite{PW11} a more direct algebraic approach, albeit a key tool \cite{HN00} used herein is essentially analytic.
\end{remark*}

\section{Preliminaries}
\vspace{.2cm}
Recall the following standard definitions from \cite{PW11} with Definition 2.2 evolving from \cite{PW11}.

\bD{D:d1}Let $J$ be an ideal of $B(H)$ (i.e., a $B(H)$-ideal) and $S \in J$.

\textbullet \, The principal $B(H)$-ideal generated by the single operator $S$ is given by
\begin{center}
$\left(S\right)$ := $\bigcap \{I \mid I$\, is\, a\, $B(H)$-ideal\, containing\, $S\,\}$
\end{center}

\textbullet \, The principal $J$-ideal generated by $S$ is given by
\begin{center}
$\left(S\right)_J$ := $\bigcap \{\mathcal I \mid \mathcal I$ is a $J$-ideal containing $S\,\}$
\end{center}

\textbullet \, As above for principal $J$-ideals, likewise for an arbitrary subset $\mathscr S \subset J$, $(\mathscr S)$ and $(\mathscr S)_J$ respectively denote the smallest $B(H)$-ideal and the smallest $J$-ideal generated by the set $\mathscr S$. In particular, $(\mathscr S) = (\mathscr S)_{B(H)}$. Denote $(\mathcal I)$ as the $B(H)$-ideal generated by the $J$-ideal $\mathcal I$.  
\eD

\bD{D:d2}
For a $J$-ideal $\mathcal I$, the algebraic $J$-interior of $\mathcal I$  is denoted by $\mathcal I^0 := \text{span}\{\mathcal IJ + J\mathcal I\} + J(\mathcal I)J$ where $\mathcal IJ := \{ AB \mid A \in \mathcal I,\,B\in J\}$, $J\mathcal I$ is defined similarly, and the ideal product $J(\mathcal I)J$ is the $B(H)$-ideal given by  $J(\mathcal I)J = \{ES'F \mid E,\,F\in J,\, S' \in (\mathcal I)\}$ (single triple operator products) where equality follows from \cite[Lemma 6.3]{DFWW}.
\eD

\bR{R:2.1}
For the $J$-ideal $\mathcal I = (\mathscr S)_J$ generated by a set $\mathscr S$, $\mathcal I^0$ has the simpler form: \begin{center}
$(\mathscr S)_J^0 = span\{\mathscr SJ + J\mathscr S\} + J(\mathscr S)J$
\end{center}
\noindent Notice also that $(\mathscr S)_J^0 \subset (\mathscr S)J$.

Indeed, for $\mathscr S = \{S_\alpha\}_{\alpha \in \mathcal A}$,  $(\mathscr S)_J = \displaystyle{\bigcup_{\{\alpha_1,\cdots,\alpha_j\} \subset \mathcal A}}\{(S_{\alpha_{1}})_J + (S_{\alpha_{2}})_J + \cdots + (S_{\alpha_{j}})_J\}$ where $\mathcal A$ is an index set. 
By definition, $(\mathscr S)_J^0 = \text{span}\{(\mathscr S)_JJ+ J(\mathscr S)_J\}+ J((\mathscr S)_J)J$, and because $((\mathscr S)_J) = (\mathscr S)$ one has the simplification: $(\mathscr S)_J^0 = \text{span}\{(\mathscr S)_JJ+ J(\mathscr S)_J\}+ J(\mathscr S)J$. 
The algebraic structure for principal $J$-ideals  yields $(S_\alpha)_J = \mathbb CS_\alpha + JS_\alpha + S_\alpha J + J(S_\alpha)J$ for each $\alpha \in \mathcal A$. 
So $\text{span}\{(S_\alpha)_JJ + J(S_\alpha)_J\} + J(S_\alpha)J\subset S_\alpha J + JS_\alpha + J(S_\alpha)J$ for each $\alpha \in \mathcal A$. 
Therefore $ (\mathscr S)_J^0 \subset span\{\mathscr SJ + J\mathscr S\} + J(\mathscr S)J$, and since the reverse inclusion is obvious, one has equality. 
\eR
We recall here the definition of $J$-softness of $B(H)$-ideals \cite[Definition 2.5]{PW11}.
\bD{D:1} \quad \\
For $B(H)$-ideals $I$ and $J$, the ideal $I$ is called ``$J$-soft" if $IJ = I$. \\
Equivalently in the language of s-numbers:\\ 
For every $A \in I$, $s_n(A) = O(s_n(B)s_n(C))$ for some $B \in I, C \in J, m \in \mathbb N$.
\eD
\noindent Because $IJ \subset J$, only $B(H)$-ideals that are contained in $J$ can be $J$-soft.
\newpage
\bR{R:1}
\emph{Standard facts on $B(H)$-ideals from \cite[Remark 2.2]{PW11}.}\\
\item (i) \cite[Sections 2.8, 4.3]{DFWW} (see also \cite[Section 4]{KW07}): If $I,J$ are $B(H)$-ideals then the product $IJ$, which is both associative and commutative, is the $B(H)$-ideal given by the characteristic set $\Sigma(IJ) = \{\xi \in c_o^* \mid \xi \leq \eta\rho$ for some $\eta \in \Sigma(I)$ and $\rho \in \Sigma(J)\}$.
In abstract rings, the ideal product is defined as the class of finite sums of products of two elements, $IJ := \{\displaystyle{\sum_{finite}} a_ib_i \mid a_i \in I, b_i \in J\}$, but in $B(H)$ the next lemma shows finite sums of operator products defining $IJ$ can be reduced to single products. \\
\item (ii)  \cite[Lemma 6.3]{DFWW} Let $I$ and $J$ be proper ideals of $B(H)$. If $A \in IJ$, then $A = XY$ for some $X \in I$ and $Y \in J$.\\
\item (iii)  \cite[Section 1]{KW07} For $T \in B(H)$, $s(T)$ denotes the sequence of s-numbers of $T$.
Then $ A \in (T)$ if and only if $s(A) = \textrm{O}(D_m(s(T))$ for some $m \in \mathbb{N}$.
Moreover, for a $B(H)$-ideal $I$, $A \in I$  if and only if $A^* \in I$ if and only if $|A| \in I$ (via the polar decomposition $A = U|A|$ and $U^*A = |A| = (A^*A)^{1/2}$), with all equivalent to $\diag s(A) \in I$.\\
\item (iv)  The lattice of $B(H)$-ideals forms a commutative semiring with multiplicative identity $B(H)$. 
That is, the lattice is commutative and associative under ideal addition and multiplication (see \cite[Section 2.8]{DFWW}) and it is distributive.
Distributivity with multiplier $K(H)$ is stated without proof in \cite[Lemma 5.6-preceding comments]{KW07}.
The general proof is simple and is as follows. For $B(H)$-ideals $I,J,K$, one has
$I,J \subset I+J := \{A+B \mid A \in I, B \in J\}$ and so $IK, JK \subset (I+J)K$, so one has $IK + JK \subset (I+J)K$.
The reverse inclusion follows more simply if one invokes (ii) above: $X \in (I+J)K$ if and only if $X = (A+B)C$ for some $A \in I, B \in J, C \in K$.
The lattice of $B(H)$-ideals is not a ring because, for instance, $\{0\}$ is the only $B(H)$-ideal with an additive inverse, namely, $\{0\}$ itself, so it is not an additive group. It is also clear that $B(H)$ is the multiplicative identity but no $B(H)$-ideal has a multiplicative inverse.\\
\eR
We summarize the main results of \cite{PW11} generalizing the 1983 work of Fong and Radjavi and characterizing all finitely generated linear $J$-ideals. 
Though not needed here, we note that \cite{PW11} provided similar results for real-linear and nonlinear $J$-ideals. 
\begin{theorem}\label{T:2.6}\cite[Theorem 4.5]{PW11}
For $\mathscr S := \{S_{1}, \cdots, S_{N}\} \subseteq J$, the following are equivalent.

(i) The finitely generated $J$-ideal $(\mathscr S)_J$ is a $B(H)$-ideal.

(ii) The $B(H)$-ideal $(\mathscr S)$ is $J$-soft, i.e., $(\mathscr S)$ = $J(\mathscr S)$ (equivalently, $(\mathscr S)$ = $(\mathscr S)J$).

(iii) For all $1 \leq j \leq N$, $S_{j} = \displaystyle{\sum_{k=1}^{N}}\displaystyle{\sum_{i=1}^{n{(j,k)}}}A_{ijk}S_{k}B_{ijk}$
 
\qquad \text{for some} $A_{ijk}, B_{ijk} \in J, n(j,k) \in \mathbb{N}$.

(iv)  For all $1 \leq j \leq N$, 
$s(S_{j}) = \textrm{O}(D_{m}(s(|S_{1}|+ \cdots +|S_{N}|))s(T)) \text{~for some~} T \in J \text{~and~} m \in \mathbb{N}$.
 \end{theorem}

 \begin{theorem}\label{T:2.7}\cite[Theorem 4.6]{PW11} (Structure theorem for finitely generated $J$-ideal $(\mathscr S)_J$ for $\mathscr S = \{S_1, \cdots, S_N\}$)\\
The algebraic structure of the finitely generated $J$-ideal $(\mathscr S)_J$ is given by 
\begin{equation*}
(\mathscr S)_J = \mathbb{C}S_{1} + \cdots + \mathbb{C}S_{N} + JS_{1} + \cdots + JS_{N} +  S_{1}J + \cdots + S_{N}J + J(\mathscr S)J
\end{equation*}
So one has
\begin{equation*}
J(\mathscr S)J \subseteq JS_{1} + \cdots + JS_{N} +  S_{1}J + \cdots + S_{N}J + J(\mathscr S)J \subseteq (\mathscr S)_J \subseteq (\mathscr S)
\end{equation*}
which first two, $J(\mathscr S)J$ and $JS_{1} + \cdots + JS_{N} +  S_{1}J + \cdots + S_{N}J + J(\mathscr S)J$ respectively, are a $B(H)$-ideal and a $J$-ideal.
The inclusions collapse to merely  
\begin{equation*}
J(\mathscr S)J = (\mathscr S)_J = (\mathscr S)
\end{equation*}
if and only if the finitely generated $B(H)$-ideal $(\mathscr S)$ is $J$-soft. 
\end{theorem}

\newpage

\section{The Hamel dimension of ideals}
In this section we show that the Hamel dimension and the cardinality of every nonzero $J$-ideal is precisely $\mathfrak c$. 
This will impact the codimension of the algebraic $J$-interior $\mathcal I^0$ in $\mathcal I$ and lead to questions on the possible generating sets for general $J$-ideals (see Question \ref{Q:11} and Section 6--Questions 1, 2, 5).

It is straightforward to see that $\mathcal I^0 := \text{span}\{\mathcal IJ + J\mathcal I\} + J(\mathcal I)J$ is an ideal of $\mathcal I$ and hence is a complex vector subspace of $\mathcal I$. 
The quotient space $\mathcal I/\mathcal I^0$ is a complex vector space and therefore has a Hamel basis where the Hamel dimension is invariant over all Hamel bases. 
The key notion used in our results is the Hamel dimension of  $\mathcal I/\mathcal I^0$ relative to its vector space structure. ($\mathcal I^0$ being an ideal of $\mathcal I$, the quotient space $\mathcal I/\mathcal I^0$ is also a ring but we will not exploit the ring structure.)
 
\bP{P:0.1}
For the $J$-ideal $(\mathscr S)_J$ generated by a set $\mathscr S$ and $(\mathscr S)_J^0 =\text{span} \{\mathscr SJ + J\mathscr S\} + J(\mathscr S)J$, the Hamel dimension of the quotient space $(\mathscr S)_J/(\mathscr S)_J^0$ is at most the cardinality of the generating set $\mathscr S$.
\eP
\bp
From general ring theory, for $\mathscr S = \{S_\alpha\}_{\alpha \in \mathcal A}$,  $(\mathscr S)_J = \displaystyle{\bigcup_{\{\alpha_1,\cdots,\alpha_j\} \subset \mathcal A}}\{(S_{\alpha_{1}})_J + (S_{\alpha_{2}})_J + \cdots + (S_{\alpha_{j}})_J\}$ where $\mathcal A$ is an index set.
Combining this with the algebraic structure for principal $J$-ideals implied by Theorem 2.7 (or for principal $J$-ideals in particular, see also \cite[Proposition 4.2]{PW11}), one obtains $(\mathscr S)_J = \text{span}\, \mathscr S+ (\mathscr S)_J^0$. 
For $(\mathscr S)_J^0$ being a linear subspace of $(\mathscr S)_J$, one can show that the quotient space $(\mathscr S)_J/(\mathscr S)_J^0 = \text{span}\{[S_\alpha]\}$ where $\alpha \in \mathcal A$, $ S_\alpha \in \mathscr S$ and $[S_\alpha]$ denotes its quotient space coset. Therefore the Hamel dimension of the vector space $(\mathscr S)_J/(\mathscr S)_J^0 $ is at most the cardinality of $\mathscr S$.
\ep

Finishing up our discussion on the Hamel dimension, the following proposition which we need in Example \ref{E:11} is probably a well-known fact but we include it here for completeness.
\bP{P:2.1}
The Hamel dimension of $F(H)$, the $B(H)$-ideal of finite rank operators, is $\mathfrak c$ when $H$ is separable (at least $\mathfrak c$ for $H$ non-separable).
\eP
\bp
Cardinal arithmetic applied to matrices  when $H$ is separable shows $|F(H)| \leq |B(H)| \leq \mathfrak c$, so the Hamel dimension of $F(H)$ is at most $\mathfrak c$.
Suppose the Hamel dimension of $F(H)$ is strictly less than $\mathfrak c$. Let $\mathcal B := \{F_\alpha \in F(H)\,| \alpha \in \mathcal A\}$ be a  Hamel basis for $F(H)$ with cardinality $|\mathcal A| < \mathfrak c$. 
Denote a finite basis for the range of $F_\alpha$ by $B_\alpha$. 
So $|\displaystyle{\bigcup_{\alpha\in \mathcal A}}B_\alpha| = |\mathcal A|$ and from the set $\displaystyle{\bigcup_{\alpha\in \mathcal A}} B_\alpha \subset H$, one can extract a maximal linearly independent set $E$ of cardinality at most $|\mathcal A|<\mathfrak c$. 
Since the Hamel dimension of  infinite-dimensional Hilbert space is at least $\mathfrak c$ \cite[Lemma 3.4]{HN00}, there is a $0\neq f \in H$  for which the set $E \bigcup \{f\}$ forms a linearly independent set. 
Consider the rank one operator $f \bigotimes f$. Since $\mathcal B$ is a Hamel basis for $F(H)$,  $f \bigotimes f = \displaystyle{\sum_{i=1}^{n}a_{i}F_{\alpha_{i}}}$ for some $a_{i} \in \mathbb C, n \in \mathbb N$. 
So, in particular, $f \bigotimes f (f) = \displaystyle{\sum_{i=1}^{n}a_{i}F_{\alpha_{i}}}(f)$ which implies $0 \neq\left<f,f\right>f = \displaystyle{\sum_{i=1}^{n}a_{i}F_{\alpha_{i}}(f)} = \displaystyle{\sum_{i=1}^{m}}b_{\beta_{i}}e_{\beta_{i}}$ where $e_{\beta_{i}} \in E \subset \displaystyle{\bigcup_{\alpha\in \mathcal A}} B_\alpha$, hence $f\in \text{span}\, E$ contradicting that $E \bigcup \{f\}$ is a linearly  independent set. 
In summary, the assumption that the Hamel dimension of $F(H)$ is strictly less than $\mathfrak c$ led to the existence of this $f$ and hence to this contradiction. 
So the Hamel dimension of $F(H)$ is precisely $\mathfrak c$, and consequently the cardinality $|B(H)| = \mathfrak c$.
\ep
An unrelated and interesting question on Hamel bases appears in \cite{B}. 
\bR{R:2}
The Hamel dimension of every nonzero $J$-ideal $\mathcal I$ is precisely $\mathfrak c$. 
Indeed, $F(H) \subset \mathcal I$ (see \cite[Section 6, \P 2]{PW11}) so  Proposition \ref{P:2.1} implies the Hamel dimension of $\mathcal I$ is at least $\mathfrak c$. 
Also, since $\mathcal I \subset B(H)$ and since cardinality of $B(H)$ is precisely $\mathfrak c$, the Hamel dimension of $\mathcal I$ is at most $\mathfrak c$. 
Hence the Hamel dimension of $\mathcal I$ is precisely $\mathfrak c$. 
Moreover, because $F(H) \subset \mathcal I \subset B(H)$, $\mathfrak c = |F(H)| \le |\mathcal I| \le |B(H)| = \mathfrak c$ so the cardinality of $\mathcal I$ is also precisely $\mathfrak c$.
\eR
\newpage
\section{Main Results: Structure and Softness}

As mentioned earlier in Section 1, one of our main contributions in \cite{PW11} was the generalization of Fong and Radjavi's result \cite[Theorem 2]{FR83} by showing that the principal $J$-ideals and the finitely generated $J$-ideals that are also $B(H)$-ideals must be $J$-soft. 
Here we show the same for $J$-ideals generated by countable sets and, absent CH, generated by sets of cardinality strictly less than $\mathfrak c$. 
We also show that if a $J$-ideal $\mathcal I$ is generated by sets of cardinality equal to $\mathfrak c$, then $\mathcal I$ being a $B(H)$-ideal does not necessarily imply that $\mathcal I$ is $J$-soft (Example \ref{E:11}). 
Our main softness theorem is:
\begin{theorem}\label{T:1}
 The $J$-ideal $(\mathscr S)_J$ generated by sets $\mathscr S$ of cardinality strictly less than $\mathfrak c$ is a $B(H)$-ideal if and only if  the $B(H)$-ideal $(\mathscr S)$ is $J$-soft. 
Moreover, for a $J$-ideal $\mathcal I$  with the Hamel dimension of $\mathcal I/\mathcal I^0$ strictly less than $\mathfrak c$, $\mathcal I$ is a $B(H)$-ideal if and only if $(\mathcal I)$ is $J$-soft, where $(\mathcal I)$ is the $B(H)$-ideal generated by $\mathcal I$.
\end{theorem}
\bp

$\Rightarrow$: Since $(\mathscr S)J \subset (\mathscr S)$, it suffices to show $(\mathscr S) \subset (\mathscr S)J$. 
Assume otherwise that there is some $T \in (\mathscr S) \setminus (\mathscr S)J$. 
We claim $(\mathscr S) = (\mathscr S)_J$ so that $T \in (\mathscr S)_J$. 
Since every $B(H)$-ideal is also a $J$-ideal, $(\mathscr S)$ is a $J$-ideal containing $\mathscr S$. 
And $(\mathscr S)_J$ being the smallest $J$-ideal containing $\mathscr S$, one has $(\mathscr S)_J \subset (\mathscr S)$. 
The minimality of $(\mathscr S)$ as a $B(H)$-ideal containing $\mathscr S$  and $(\mathscr S)_J$ assumed to be a $B(H)$-ideal imply $(\mathscr S) \subset (\mathscr S)_J$. 
Hence $(\mathscr S)_J = (\mathscr S)$ and therefore $T \in (\mathscr S)_J$ or, equivalently because when $(\mathscr S)_J = (\mathscr S)$ it is a $B(H)$-ideal, one has the equivalent condition $\diag{s(T)} \in (\mathscr S)_J$.

Using the s-number sequence $s(T)$ we now construct a sequence of operators $D_n$ as follows. 
For each $n \geq 1$, let $D_n$  be the diagonal operator with $s_{2^{n-1}(2k-1)}(T)$ at the $2^{n-1}(2k-1)$ scattered diagonal positions for $k \geq 1$ and with  zeros elsewhere. 
Every positive integer has this unique product decomposition $2^{n-1}(2k-1)$. 
Notice then that the diagonal sequences of the $D_n$'s have  pairwise disjoint support and the formal direct sum $\sum^\oplus D_n$ (which incidentally converges in the operator norm) is precisely  $\diag {s(T)}$. 
Recall that $(\mathscr S)^{0}_J$ is a $J$-ideal and a complex vector subspace of $(\mathscr S)_J$ so their quotient $(\mathscr S)_J/(\mathscr S)_J^0 $ is a vector space. 
We will use these $D_n$'s to imbed isomorphic copies of  $\ell^p$ (for every $1\leq p \leq \infty$) inside the quotient space  $(\mathscr S)_J/(\mathscr S)_J^0 $. (In fact, we imbed isometric isomorphic copies of $\ell^p$ inside the quotient space  $(\mathscr S)_J/(\mathscr S)_J^0 $.)

We next show that for each $n \geq 1$, $(D_n) = (T)$ and that $D_n \notin (\mathscr S)^{0}_J$.
Clearly $(D_n) \subset (T)$ since $D_n = PT$ for a suitable diagonal projection operator, so it remains to show $(T) \subset (D_n)$.
For each $n\geq 1$, $D_n$ is explicitly given by $\diag{(\underbrace{0, \cdots, 0}_{(2^{n-1}-1)- times}, s_{2^{n-1}}(T), \underbrace{0, \cdots, 0}_{(2^{n}-1)-times}, s_{3\cdot2^{n-1}}(T),\underbrace{0, \cdots, 0}_{(2^{n}-1)-times},s_{5\cdot2^{n-1}}(T),\underbrace{0, \cdots, 0}_{(2^{n}-1)-times}, \cdots)}$\\
so the $2^n$-fold ampliation of $s(D_n)$ is 
 $(\underbrace{s_{2^{n-1}}(T), \cdots,s_{2^{n-1}}(T)}_{2^n-times},\underbrace{s_{3\cdot2^{n-1}}(T), \cdots,s_{3\cdot2^{n-1}}(T)}_{2^n-times}, \cdots) \in \Sigma((D_n))$ \\and since 
$(s_{2^{n-1}+1}(T),s_{2^{n-1}+2}(T), \cdots) \leq (\underbrace{s_{2^{n-1}}(T), \cdots,s_{2^{n-1}}(T)}_{2^n-times},\underbrace{s_{3\cdot2^{n-1}}(T), \cdots,s_{3\cdot2^{n-1}}(T)}_{2^n-times}, \cdots)$,\\
so $\left<s_{2^{n-1}+k}(T)\right>_{k=1}^{\infty}~ \in \Sigma((D_n))$ and hence $(\underbrace{0, \cdots, 0}_{(2^{n-1})-times},s_{2^{n-1}+1}(T), s_{2^{n-1}+2}(T), \cdots) \in \Sigma((D_n))$.\\ 
But also the finitely supported sequence $(s_1(T), s_2(T), \cdots, s_{2^{n-1}}(T), 0, \cdots ) \in \Sigma((D_n))$. Adding both sequences one obtains precisely $s(T)$. 
Then since $\Sigma((D_n))$ is additive, we have $s(T) \in \Sigma((D_n))$  and hence $\diag{s(T)} \in (D_n)$. 
Therefore $(T)\subset (D_n)$, and then from the reverse inclusion above one has $(T) = (D_n)$. To see that $D_n \notin (\mathscr S)^{0}_J$, assume otherwise that $D_n \in (\mathscr S)^{0}_J$. 
Then $D_n \in (\mathscr S)^{0}_J \subset (\mathscr S)J$ so $T \in (D_n) \subset(\mathscr S)J$ contradicting the assumption $T \in (\mathscr S) \setminus  (\mathscr S)J$. 
So $D_n \notin (\mathscr S)^{0}_J$ for all $n\geq1$. (Other choices of the diagonal sequences for the $D_n$'s are possible. 
Besides disjoint supports or ``almost'' disjoint supports, the only feature needed is bounded gaps between their nonzero entries.)

The set $X_p := \{\sum^\oplus a_iD_i\mid ||\left<a_i\right>||_p < \infty,\, a_i \in \mathbb C\} \subset (\mathscr S)_J$ for $1\leq p \leq \infty$. 
The  inclusion is because the s-number sequence $s(\sum^\oplus a_iD_i) \leq  ||\left<a_i\right>||_\infty \,s(T)$ and because $(\mathscr S)_J$ contains $\diag{s(T)}$ and $(\mathscr S)_J$, being assumed a $B(H)$-ideal, is hereditary. 
Clearly $X_p$, with its cannonical  vector space structure, is a linear subspace of $(\mathscr S)_J$.
Under the natural projection map, the set of cosets of elements of $X_p$ in the quotient space $(\mathscr S)_J/(\mathscr S)_J^0$ is given by the linear subspace $X_p^{'} := \{[\sum^\oplus a_iD_i]\mid ||\left<a_i\right>||_p <\infty \}$. 

We first show that the map $\sum^\oplus a_iD_i \rightarrow [\sum^\oplus a_iD_i]$ is a one-to-one map, that is, each coset $[\sum^\oplus a_iD_i]$ has a unique element of the form $\sum^\oplus a_iD_i$. 
Indeed, if $[\sum^\oplus a_iD_i]=[\sum^\oplus a_i^{'}D_i]$, then $[\sum^\oplus(a_i - a_i^{'})D_i] = [0]$. 
This implies $\sum^\oplus(a_i - a_i^{'})D_i \in (\mathscr S)_J^0 \subset (\mathscr S)J$. 
Suppose there exist $i_0$ such that $a_{i_0} \neq a_{i_0}^{'}$. 
Since $(\mathscr S)J$ is a $B(H)$-ideal, multiplying $\sum^\oplus(a_i - a_i^{'})D_i$ by a suitable projection it follows that the diagonal operator $(a_{i_0} - a_{i_0}^{'})D_{i_0} \in (\mathscr S)J$, and hence $D_{i_0} \in (\mathscr S)J$. 
But $(T) = (D_{i_0})\subset (\mathscr S)J$ implying $T \in (\mathscr S)J$, again contradicting $T \in (\mathscr S)\setminus(\mathscr S)J$. 
Therefore $a_i = a_i^{'}$ for all $i \geq 1$. 
So the map $\sum^\oplus a_iD_i \rightarrow [\sum^\oplus a_iD_i]$ is a one-to-one map which is clearly linear, 
and therefore the map $\left<a_i\right> ~\rightarrow [\sum^\oplus a_iD_i]$ is an isomorphism from $\ell^p$ onto $X_p^{'}$. 
(In fact, it is straightforward to show that $||[\sum^\oplus a_iD_i]||:= ||\left<a_i\right>||_p$ is a well-defined complete norm on $X_p^{'}$ which establishes that this linear map is an isometric isomorphism under this induced norm on $X_p^{'}$, but we will not exploit this isometric property.)
 
$\ell^p$ being an infinite-dimensional Banach space over the complex numbers, the cardinality of a Hamel basis for $\ell^p$ is at least $\mathfrak c$ (\cite[Lemma 3.4]{HN00}).
Then likewise a Hamel basis $\mathcal B^{'}$ for $X_p^{'}$ is at least $\mathfrak c$, since isomorphisms preserve Hamel bases.
Since $X_p^{'}$ is a vector subspace of  $(\mathscr S)_J/(\mathscr S)_J^0 $, every Hamel basis of a subspace can be extended to a Hamel basis of the full space and because the cardinality of all Hamel bases of a vector space is invariant, it follows that $|\mathcal B^{'}| \leq |\mathcal B|$ for $\mathcal B$ a Hamel basis of  $(\mathscr S)_J/(\mathscr S)_J^0 $. 
Also since the generating set $\mathscr S$ for $(\mathscr S)_J$ has cardinality strictly less than $\mathfrak c$, $|\mathcal B| < \mathfrak c$ by Proposition \ref{P:0.1}.
Therefore $\mathfrak c \leq |\mathcal B^{'}| \leq |\mathcal B| < \mathfrak c$, a set theoretic contradiction. 
To sum up, this contradiction followed from assuming properness of the inclusion $(\mathscr S)J \subsetneq (\mathscr S)$. 
Therefore $(\mathscr S)J = (\mathscr S)$, that is, $(\mathscr S)$ is $J$-soft. 

Next we prove the first implication of the second assertion of this theorem, that is, if $\mathcal I$ is a $B(H)$-ideal, then $(\mathcal I)$ is $J$-soft. Following the same method as used above for $\mathcal I=(\mathscr S)_J$, notice that the contradiction arose from assuming properness of the inclusion $(\mathcal I)J \subsetneq (\mathcal I)$, that is, we showed there how the assumption of  $(\mathcal I)J \neq(\mathcal I)$ led to an imbedding of $X_{p}^{'}$ (an isometric isomorphic copy of $\ell^p$) into $\frac{\mathcal I}{\mathcal I^0}$ without depending on cardinality of $\mathscr S$, and yet still violating  dim $\frac{\mathcal I}{\mathcal I^0} < \mathfrak c$. 

$\Leftarrow$: From general ring theory, for $\mathscr S = \{S_\alpha\}_{\alpha \in \mathcal A}$:~$(\mathscr S)_J = \displaystyle{\bigcup_{\{\alpha_1,\cdots,\alpha_j\} \subset \mathcal A}}\{(S_{\alpha_{1}})_J + (S_{\alpha_{2}})_J + \cdots + (S_{\alpha_{j}})_J\}$\\ where $\mathcal A$ is an index set with $|\mathcal A|< \mathfrak c$. Using the algebraic structure for principal $J$-ideals implied by Theorem 2.7 (or for principal $J$-ideals in particular, see also \cite[Proposition 4.2]{PW11}), one obtains $(\mathscr S)_J = \text{span}\,\mathscr S + (\mathscr S)_J^0$. Using Remark \ref{R:2.1} for $(\mathscr S)_J^0$, one obtains $(\mathscr S)_J = \text{span}\,\mathscr S + \text{span}\{\mathscr SJ + J\mathscr S\} + J(\mathscr S)J$. 
Then $(\mathscr S)\subset(\mathscr S)_J$ because $(\mathscr S) = (\mathscr S)J$ and the commutativity of $B(H)$-ideal multiplication implies $(\mathscr S) = J(\mathscr S)J$ \cite[Sections 2.8, 4.3]{DFWW}. 
But the minimality of $(\mathscr S)_J$  as a $J$-ideal containing $\mathscr S$ implies $(\mathscr S)_J \subset (\mathscr S)$, and so combining both inclusions, one obtains $(\mathscr S) = (\mathscr S)_J$, that is, $(\mathscr S)_J$ is a $B(H)$-ideal. 

Finally, we prove the second implication of the second assertion of this theorem, that is, if $(\mathcal I)$ is $J$-soft, then $\mathcal I$ is a $B(H)$-ideal. Notice that when $(\mathcal I)$ is $J$-soft, one obtains $(\mathcal I) = (\mathcal I)J = J(\mathcal I)J$, hence  $\mathcal I \subset J(\mathcal I)J$. And $J(\mathcal I)J \subset \mathcal I$ because $\mathcal I$ is a $J$-ideal. Combining both inclusions one obtains $\mathcal I = J(\mathcal I)J$ which is a $B(H)$-ideal.  
\ep

\bR{R:4}In Theorem \ref{T:1} only one implication (i.e., $(\mathscr S)_J$ is a $B(H)$-ideal $\Rightarrow$ $(\mathscr S)$ is $J$-soft) requires the cardinality of $\mathscr S$ to be strictly less than the continuum. The reverse implication holds for arbitrary $\mathscr S$. 
\eR
 \bQ{Q:11} Is the codimension condition equivalent to the $J$-ideal $\mathcal I$ possessing a generating set of cardinality less than $\mathfrak c$. (See also Section 6, Question 5.)
\eQ

As a consequence of Theorem \ref{T:1} we have

 \bT{T:3}(Structure theorem for  $(\mathscr S)_J$ when $|\mathscr S| < \mathfrak c$)\\
The algebraic structure of the $J$-ideal $(\mathscr S)_J$ generated by the set $\mathscr S$ is given by
\begin{center}
$(\mathscr S)_J = \text{span}\{\mathscr S + J\mathscr S + \mathscr SJ\} + J(\mathscr S)J$
\end{center}
$J(\mathscr S)J $ is a $B(H)$-ideal, $\text{span}\{J\mathscr S + \mathscr SJ\} + J(\mathscr S)J$ is a $J$-ideal, and $J(\mathscr S)J \subset \text{span}\{J\mathscr S + \mathscr SJ\} + J(\mathscr S)J \subset (\mathscr S)_J$.
This inclusion collapse to $J(\mathscr S)J = (\mathscr S)_J$ if and only if  $(\mathscr S)$ is $J$-soft.
\eT
\vspace{.2cm}
The fact that the cardinality of every nonzero $J$-ideal $\mathcal I$ is $\mathfrak c$ (Remark \ref{R:2}) implies that generating sets for  $\mathcal I$ have at most $\mathfrak c$ elements. 
So in view of Theorem \ref{T:1}, the only $J$-softness cases left to investigate are: if $\mathcal I$ cannot be generated by fewer than $\mathfrak c$ elements or (possibly more general, see Question \ref{Q:11}) at least if the Hamel dimension of the quotient space $\mathcal I/\mathcal I^0$ is equal to $\mathfrak c$, does either of these imply $\mathcal I$ is $J$-soft? 
Indeed, we show in Example \ref{E:11} that Theorem \ref{T:1} is the best possible result of its type by giving an example of a $J$-ideal that is also a $B(H)$-ideal which is not $J$-soft. 
By the contrapositive of Theorem \ref{T:1}, this $J$-ideal has no generating sets of cardinality less than $\mathfrak c$ and the Hamel dimension of its quotient $\mathcal I/\mathcal I^0$ is precisely $\mathfrak c$.
 
\bE{E:11}
Consider $J = K(H)$ and $\mathcal I = (\diag{\left<\frac{1}{n}\right>})$ the principal $B(H)$-ideal generated by the diagonal operator $\diag{\left<\frac{1}{n}\right>}$. 
Since every $B(H)$-ideal is a $J$-ideal, $\mathcal I$ is also a  $J$-ideal. 
We will show that the Hamel dimension of the quotient space $\mathcal I/\mathcal I^0$ is precisely $\mathfrak c$, but yet $\mathcal I$ is not $K(H)$-soft. 
Indeed, $\mathcal I^0 = \text{span}\{\mathcal IK(H) + K(H)\mathcal I\} + K(H)\mathcal IK(H)$ and one can show that $\diag{\left<\frac{1}{n}\right>} \notin \mathcal I^0$ \cite[Example 3.3]{PW11}. 
In the proof of Theorem \ref{T:1} taking $T = \diag {\left<\frac{1}{n}\right>}$, imbed $X_p'$ into $\mathcal I/\mathcal I^0$ for any $1\leq p \leq \infty$. 
So the Hamel dimension of the quotient space $\mathcal I/\mathcal I^0$ is at least $\mathfrak c$. Since $\mathcal I$ is a nonzero $J$-ideal, the Hamel dimension of $\mathcal I$ is equal to $\mathfrak c$ (Remark \ref{R:2}), so the Hamel dimension of the quotient space $\mathcal I/\mathcal I^0$ is at most $\mathfrak c$. 
 Therefore the Hamel dimension of  $\mathcal I/\mathcal I^0$ is equal to $\mathfrak c$. 
But we know that $(\diag{\left<\frac{1}{n}\right>})$ is not $K(H)$-soft \cite[Example 3.3]{PW11}, but for completeness we repeat the proof here.
If it were $K(H)$-soft, then $(\diag\left<\frac{1}{n}\right>) = (\diag\left<\frac{1}{n}\right>)K(H)$
which further implies\\ $\left<\frac{1}{n}\right> \in \Sigma ((\diag\left<\frac{1}{n}\right>)K(H))$, i.e.,  
$\left<\frac{1}{n}\right> = o(D_{m}\left<\frac{1}{n}\right>)$ for some $m \in \mathbb N$, contradicting
$\left(\frac{\left<\frac{1}{n}\right>}{D_{m}\left<\frac{1}{n}\right>}\right)_k = \frac{\frac{1}{mj+r}}{\frac{1}{j+1}} \rightarrow \frac{1}{m}$
as $k \rightarrow \infty$ where $k = mj+r$.
 
\eE
 
\section{Questions and results on $J$-ideals}
\vspace{.2cm}
In this section we address some of the questions posed in \cite[Sections 6-7]{PW11} and pose new questions.\\

The algebraic structure of a principal $J$-ideal generated by $S$ is $(S)_J = \mathbb CS + SJ + JS + J(S)J$. For idempotent $B(H)$-ideals $J$ (i.e., $J^{2} = J$),  $J(S)J = SJ + JS + J(S)J$ since $J(S)J = (S)J^2 = (S)J$ and $SJ, JS \subset (S)J = J(S)$. 
So the algebraic structure of $(S)_J$ simplifies to $(S)_J = \mathbb CS + J(S)J$.  
Is it necessary for $J$ to be idempotent for $J(S)J = SJ + JS + J(S)J$ to hold? (This is related to \cite[Remark 6.4, Section 7-Question 6]{PW11}: Find a necessary and sufficient condition(s) to make $J(S)J = JS + SJ + J(S)J$.) 
The following example shows that $J$ need not be idempotent for the equality to hold.

\bE{E:1}
For $0< p < \infty$, the Schatten p-ideal $C_p$ is not idempotent because $C_p^2 = C_{p/2}\neq C_p$. Moreover, for $S = \diag {\left<\frac{1}{n^n}\right>}\in C_p$ we claim that $C_p(S)C_p = SC_p + C_pS + C_p(S)C_p$ which will follow below from the $C_p$-softness of $(S)$. 
Indeed, for $m > 2$, recall that $(D_m(s(S)))_n = s_{\lceil \frac{n}{m}\rceil}(S)$, so 
\begin{center}
$\frac{s_n(S)}{(D_m(s(S)))_n} = \frac{(j+1)^{j+1}}{(mj+r)^{mj+r}} \leq \frac{(2j)^{2j}}{(mj)^{mj}} < \frac{(mj)^{2j}}{(mj)^{mj}} = \frac{1}{m^{(m-2)j}j^{(m-2)j}} \leq \frac{1}{j^j}$,\end{center}
 where $n = mj + r$ and the roof function $\lceil \frac{n}{m}\rceil = j$.
By the hereditary property of $\Sigma (C_p),  \frac{s(S)}{D_m(s(S))} \in \Sigma (C_p)$ which further implies that $s(S) = \frac{s(S)}{D_m(s(S))}D_m(s(S)) \in \Sigma (C_p(S))$. 
Hence $(S) = C_p(S)$ which is $C_p$-softness of $(S)$, so $(S) = C_p(S)C_p$. Because $C_p^2 \subset C_p$, it follows that $C_pS \subset C_p^2(S)C_p \subset C_p(S)C_p$ and similarly $SC_p \subset C_p(S)C_p$. 
Therefore $SC_p + C_pS \subset C_p(S)C_p$, hence $C_p(S)C_p = SC_p + C_pS + C_p(S)C_p$. 
\eE

\vspace{.5cm}
 Finishing this discussion on $JS + SJ + J(S)J$, recall that in the case of a principal $J$-ideal $(S)_J$, $JS + SJ + J(S)J$ is always a maximal $J$-ideal in $(S)_J$ \cite[Remark 6.3]{PW11}. 
The following example gives a partial answer to \cite[Section 7, Question 3]{PW11}: is  $JS + SJ + J(S)J$  always a principal $J$-ideal or is it always a non-principal $J$-ideal? 

\bE{E:2}
For $J = K(H)$ and $S = \diag {\left<\frac{1}{n}\right>}$,
\begin{center}
$JS + SJ + J(S)J = K(H)\diag {\left<\frac{1}{n}\right>} + \diag {\left<\frac{1}{n}\right>K(H)} + K(H)(\diag {\left<\frac{1}{n}\right>)}K(H)$
\end{center}
We  claim that this is not a principal $K(H)$-ideal.  
 
Suppose $JS + SJ + J(S)J = (T)_{K(H)}$. 
Since $B(H)$-ideals commute and $K(H)^2 = K(H)$, one has $K(H)(\diag {\left<\frac{1}{n}\right>)}K(H) = (\diag {\left<\frac{1}{n}\right>})K(H)$ and $K(H)\diag {\left<\frac{1}{n}\right>} + \diag {\left<\frac{1}{n}\right>K(H)} \subset (\diag {\left<\frac{1}{n}\right>})K(H)$. 
Therefore $(T)_{K(H)} = (\text{diag} \left<\frac{1}{n}\right>)K(H)$ hence $(T)_{K(H)}$ is the $B(H)$-ideal $(T)$ \cite[Theorem 2.6, for principal $J$-ideal see also Theorem 1.2]{PW11}. 
Since $T \in (\diag {\left<\frac{1}{n}\right>})K(H)$, $s(T) = \textrm {O}(D_m\left<\frac{1}{n}\right>\rho)$ for some $\rho = \left<\rho_n\right> \in c_0^*$ \cite[Section 1, \P 1]{KW07}, and without loss of generality one can assume $\rho_1 = 1$. 
Observe that the sequence $\displaystyle{\sum_{l=1}^{\infty}}\frac{1}{2^l}(D_l(\rho ^{{1/l}})) \in \Sigma (K(H)) = c_0^*$ (the cone of positive sequences decreasing to 0) because it is the $c_0$-norm limit of the  sequences $\displaystyle{\sum_{l=1}^{n}}\frac{1}{2^l}(D_l(\rho ^{{1/l}}))$ as $n \rightarrow \infty$ and $c_0$ (the sequence space of complex numbers tending to 0) is norm-closed.
So in particular, $\displaystyle{\sum_{l=1}^{\infty}}\frac{1}{2^l}(D_l(\rho ^{{1/l}}))\cdot \left<\frac{1}{n}\right> \in \Sigma(K(H)(\text{diag} \left<\frac{1}{n}\right>)) = \Sigma((T))$. 
So $\displaystyle{\sum_{l=1}^{\infty}}\frac{1}{2^l}(D_l(\rho ^{{1/l}}))\cdot {\left<\frac{1}{n}\right>} \in \Sigma((T))$, and expressing this in terms of s-numbers, for some $k, m \in \mathbb N$ and all $j = kmi + r$ for some $i$ and some $0\le r <km$,
\begin{align*}
\displaystyle{\sum_{l=1}^{\infty}}\frac{1}{2^l}(D_l(\rho ^{{1/l}}))\cdot
{\left<\frac{1}{n}\right>} &= \textrm {O}(D_ks(T))
 = \textrm {O}(D_k(D_m\left<\frac{1}{n}\right>\rho))
 = \textrm {O}(D_{km}\left<\frac{1}{n}\right>D_k(\rho)),
\end{align*}
thereby contradicting, after setting arbitrary $j = kmi +r$, for some $i$ and $0\leq r\leq km$,\\ 
$\left(\frac{\left<\frac{1}{n}\right>\displaystyle{\sum_{l=1}^{\infty}}\frac{1}{2^l}(D_l(\rho ^{1/l}))}{D_{km}\left<\frac{1}{n}\right> D_k(\rho)}\right)_j = \frac{\frac{1}{j}\displaystyle{\sum_{l=1}^{\infty}}\frac{1}{2^l}(D_l(\rho ^{{1/l}}))_j}{(D_{km}\left<\frac{1}{n}\right> D_k(\rho))_j} \geq \frac{\frac{1}{j}\displaystyle{\sum_{l=1}^{\infty}}\frac{1}{2^l}(D_l(\rho ^{{1/l}}))_j}{(D_{km}\left<\frac{1}{n}\right> D_{km}(\rho))_j}\geq \frac{\frac{1}{j}\frac{1}{2^{km}}(D_{km}(\rho ^{\frac{1}{km}}))_j}{(D_{km}(\frac{1}{n}) D_{km}(\rho))_j} =  \frac{\frac{1}{(kmi + r)} }{2^{km}\frac{1}{(i+1)} \rho_{i+1}^{1-\frac{1}{km}}}\\ = \frac{(i+1)}{2^{km}(kmi+r)\rho_{i+1}^{1-\frac{1}{km}}} $ which diverges to $\infty$ as $j \rightarrow \infty$ since $km >1$ and $\rho_i \rightarrow 0$.
\eE

\vspace{.4cm}
For principal $B(H)$-ideals $(S)$, $(T)$, $(S) = (T)$ if and only if $s(S)= \textrm{O}(D_ms(T))$ and $s(T) = \textrm{O}(D_ks(S))$ for some $m,k \in \mathbb N$. When  $s(S)$ or $s(T)$ satisfies the $\triangle_{1/2}$ condition, a simpler condition is: $(S) = (T)$ if and only if $s(S)= \textrm{O}(s(T))$ and $s(T) = \textrm{O}(s(S))$. (See \cite[Section 2, \P 1]{KW11}.) 
The following proposition gives a necessary and sufficient condition for two principal $J$-ideals to be equal \cite[Section 7, Question 4]{PW11}. We hope for a simpler condition. 
 
\bP{P:1}
For $S, T  \in J$, \begin{center}
$(S)_J = (T)_J$ \qquad if and only if \qquad  $aS + bT \in \{SJ + JS + J(S)J\} \bigcap \{TJ + JT + J(T)J\}$
\end{center}
for some nonzero $a, b \in \mathbb C$.  
\eP
\bp
$\Rightarrow$: Based on the algebraic structure of principal $J$-ideals implied by Theorem \ref{T:2.7} (or for principal $J$-ideals in particular, see also \cite[ Proposition 2.3]{PW11}), $(S)_J = (T)_J$ if and only if 
\begin{align*}
S &= \alpha T + AT + TB + \displaystyle{\sum_{i=1}^n}A_iTB_i \qquad \text{and}\qquad T = \beta S + CS + SD + \displaystyle{\sum_{j=1}^m}C_jSD_j
 \end{align*}
for some $\alpha, \beta \in \mathbb C, A, B, C, D, A_i, B_i, C_i, D_i \in J$.
If $\alpha = 0$ and $\beta \neq 0$, then substituting $S$ from the first equation in the second equation yields $T \in (T)J$ and then substituting $T$ from the second equation in the first equation yields $S \in (S)J$. 
Therefore in this case $(T) = J(T)J = (T)_J$ and $(S) = J(S)J = (S)_J$ implied by Theorem \ref{T:2.6} (or for principal $J$-ideals in particular, see also \cite[Theorem 1.2]{PW11}). 
Since $J(S)J = (S)_J = (T)_J = J(T)J$, one has $S,T \in J(S)J \bigcap J(T)J$, and because $J(S)J \bigcap J(T)J$ is a $J$-ideal, $aS + bT \in \{SJ + JS + J(S)J\} \bigcap \{TJ + JT + J(T)J\}$ for all $a, b \in \mathbb C$.  
The cases $\alpha \neq 0$ and $\beta = 0$ and $\alpha = 0 = \beta$ are handled similarly. 
Assume finally that $\alpha, \beta \neq 0$. Substituting $T$ from  the second equation in the first equation one obtains $S = \alpha\beta S + X$, where $X \in (S)J$. 
If $\alpha\beta \neq 1$, then $S \in \{SJ + JS + J(S)J\}$, and then, substituting $S$ from the first equation in  the second equation, one obtains $T = \alpha\beta T + Y$ where $Y \in (T)J$, so $T \in \{TJ + JT + J(T)J\}$. This implies $J$-softness of $(S)$ and $(T)$, hence $J(S)J = (S)_J = (T)_J = J(T)J$ and then as above $aS + bT \in \{SJ + JS + J(S)J\} \bigcap \{TJ + JT + J(T)J\}$ for all $a, b \in \mathbb C$.  
When $\alpha\beta = 1$, $\beta = \frac{1}{\alpha}$, so substituting $\beta$ in the second equation yields  
$S - \alpha T = AT + TB + \displaystyle{\sum_{i=1}^n}A_iTB_i$ and 
$S - \alpha T = -\alpha(CS + SD + \displaystyle{\sum_{j=1}^m}C_jSD_j)$.
Therefore $S - \alpha T \in \{SJ + JS + J(S)J\} \bigcap \{TJ + JT + J(T)J\}$ which is the required condition.

$\Leftarrow:$ Suppose $aS + bT \in \{SJ + JS + J(S)J\} \bigcap \{TJ + JT + J(T)J\}$ for some nonzero $a, b \in \mathbb C$. Then 
\begin{align*}
aS + bT &= AT + TB + \displaystyle{\sum_{i=1}^n}A_iTB_i \qquad \text{and}\qquad aS + bT =  CS + SD + \displaystyle{\sum_{j=1}^m}C_jSD_j
\end{align*}
for some $\alpha, \beta \in \mathbb C, A, B, C, D, A_i, B_i, C_i, D_i \in J$.
From these equalities its clear that $(S)_J = (T)_J$.
\ep
\newpage
Unlike $B(H)$-ideals, $J$-ideals do not necessarily commute as given in the following example.
\bE{E:1}
Consider $J = K(H)$ and with respect to the standard basis 
take $S$ to be the diagonal matrix 
$S$ := $\text{diag}\,(1, 0, 1/2, 0, 1/3, 0, ...)$ 
and $T$ to be the weighted shift with this same weight sequence.
We claim $(S)_{K(H)}(T)_{K(H)} \neq (T)_{K(H)}(S)_{K(H)}$.
Indeed, suppose $(S)_{K(H)}(T)_{K(H)} = (T)_{K(H)}(S)_{K(H)}$. 
Then $TS - ST \in (S)K(H)(T) = K(H)(\diag \left<\frac{1}{n^2}\right>)$. Since $TS = T$ and $ST = 0$, one has $T \in K(H)(\diag \left<\frac{1}{n^2}\right>)$.
But $T \in K(H)(\diag \left<\frac{1}{n^2}\right>)$ if and only if $\diag \left<\frac{1}{n}\right>\in K(H)(\diag \left<\frac{1}{n^2}\right>) \subset  K(H)(\diag \left<\frac{1}{n}\right>)$, the latter inclusion contradicts $\text{diag} \left<\frac{1}{n}\right> \notin K(H)(\diag \left<\frac{1}{n}\right>)$ (Example \ref{E:11} or \cite[Example 3.3]{PW11}), so $T \notin K(H)(\diag \left<\frac{1}{n^2}\right>)$, a contradiction. 
Therefore $(S)_{K(H)}(T)_{K(H)} \neq (T)_{K(H)}(S)_{K(H)}$. 
\eE
\vspace{.2cm}
\begin{center}
{QUESTIONS}\\
\end{center}
 \vspace{.2cm}
\noindent 1. Determine which $J$-ideals have generating sets of cardinality less than $\mathfrak c$?\\
 
\noindent 2. Do $J$-ideals have minimal cardinality generating sets?\\
 
\noindent 3. Traces. $B(H)$-ideals are important in part because of the importance of their traces, unitarily invariant functionals on $B(H)$-ideals. Because $B(H)$-ideals $J$ contain no unitaries this concept does not apply. Is there a modification for which one obtains a useful concept of traces for $J$-ideals?\\
 
\noindent 4. The study of commutators in $B(H)$-ideals is directly related to their traces. What can be said about the commutator structure of the $J$-ideals that can be used to motivate notions of traces for $J$-ideals?\\
 
\noindent 5. Is the dimension of $\mathcal I / \mathcal I^0$ the cardinality of some generating set for $\mathcal I$?\\
\noindent At least if dimension $\mathcal I / \mathcal I^0 < \mathfrak c $, must $\mathcal I$ possess a generating set of cardinality less than $\mathfrak c$?\\ (See also Question \ref{Q:11}.)

\end{document}